\documentclass[12pt,twoside]{article}

\usepackage{amsmath,amsfonts,amsthm,amssymb,amsopn,amstext,amscd}
\usepackage{mathabx}
\usepackage{enumerate}
\usepackage{xcolor}
\usepackage{longtable}
\usepackage{graphicx}
\usepackage{float}
\usepackage{setspace}
\usepackage{enumitem}
\usepackage{xcolor}
\usepackage{hyperref}
\usepackage{dsfont}
\usepackage[normalem]{ulem}
\allowdisplaybreaks

\newcommand{\N}{\mathbb{N}}

\newcommand{\dI}{\mathcal{I}_\Delta}
\newcommand{\I}{\mathcal{I}}
\newcommand{\T}{\mathcal{T}}

\newcommand{\cuad}{\begin{flushright}\vspace{-2ex}$\Box$\vspace{-2ex}\end{flushright}}

\newenvironment{Prf}[1][\unskip]{%
\par
\noindent
{\textbf{Proof of #1}}\newline
\vspace{-2ex}\noindent{}\newline}\cuad

\newcounter{construction}
\newenvironment{construction}[1][\unskip]{%
\refstepcounter{construction}\par
\noindent
{\textbf{Construction~\theconstruction: #1}} }

\newcommand*\samethanks[1][\value{footnote}]{\footnotemark[#1]}

\newtheorem{thm}{Theorem}
\newtheorem{prop}{Proposition}
\newtheorem{cor}{Corollary}
\newtheorem{lem}{Lemma}
\newtheorem{rem}{Remark}
\newtheorem{ex}{Examples}

\begin{document}
\pagenumbering{arabic}

\title{Defective Galton-Watson processes in a varying environment}
\author{G\"otz Kersting\thanks{This is a plain version of the paper accepted in \emph{Bernoulli} for publication (see \url{http://www.bernoulli-society.org/index.php/publications/bernoulli-journal/bernoulli-journal}). Both authors contributed equally to this work.} \thanks{Institut f\"ur Mathematik, Goethe Universit\"at, Frankfurt am Main, Germany. E-mail address: \url{kersting@math.uni-frankfurt.de}.}  \and Carmen Minuesa\samethanks[1]\  \thanks{Department of Mathematics, Autonomous University of Madrid, Madrid, Spain. E-mail address: \url{carmen.minuesa@uam.es}. }\  \thanks{Department of Mathematics, University of Extremadura, Badajoz, Spain. E-mail address: \url{cminuesaa@unex.es}. ORCID: 0000-0002-8858-3145.} }
\maketitle

\begin{abstract}
We study an extension of the so-called defective Galton-Watson processes obtained by allowing the offspring distribution to change over the generations. Thus, in these processes, the individuals reproduce independently of the others and in accordance to some possibly defective offspring distribution depending on the generation. Moreover, the defect $1-f_n(1)$ of the offspring distribution at generation $n$ represents the probability that the process hits an absorbing state $\Delta$ at that generation. We focus on the asymptotic behaviour of these processes. We establish the almost sure convergence of the process to a random variable with values in $\N_0\cup\{\Delta\}$ and we provide two characterisations of the duality extinction-absorption at $\Delta$. We also state some results on the absorption time and the properties of the process conditioned upon its non-absorption, some of which require us to introduce the notion of defective branching trees in varying environment.
\end{abstract}



\maketitle


\section{Introduction}\label{sec:intro}

A special family of branching processes, known as defective Galton-Watson processes (DGWPs), was studied in \cite{Sagitov-Minuesa-2017}. The definition of these processes is similar to classic Galton-Watson processes (GWPs), with the difference of considering defective offspring distributions. In this paper, we generalize the notion of defective Galton-Watson process by letting the defective offspring distribution change along the generations; the resulting process is called defective Galton-Watson processes  in a varying environment (DGWPVE). The dynamics of the populations described by these processes is as follows. Under the common assumption of the independence in the reproduction, each individual at generation $n-1$ produces $k$ offspring with probability $f_n[k]$, $n\in\N$, $k\in\N_0$. Moreover, the probability generating function (p.g.f.) of the offspring at generation $n$, $f_n$, satisfies $f_n(1)=\sum_{k=0}^\infty f_n[k] \leq 1$. Thus, these processes not only extend the family of DGWPs, but also the family of branching processes in a varying environment (BPVEs).


Similarly to BPVEs (see \cite{Kersting-2020}), while the presence of the varying environment is a natural assumption for modelling practical situations, it entails an added difficulty in the study of these processes. In this framework, the increase in the complexity is more noticeable since the p.g.f. of the process at the $n$-th generation, $f_{0,n}=f_1\circ \cdots\circ f_n$, does not longer satisfy $f_{0,n}(1)=1$.  
Due to all the aforementioned issues, the study of DGWPVE is challenging. Nevertheless, they constitute appropriate models for the description biological systems where each individual can show a certain trait (such as a physical feature, a mutation, or even suffering from some disease) and we are only interested in how the system evolves until the first individual shows this feature. Thus, we introduce the absorbing state $\Delta$ to denote the presence of the trait in some individual of the population and the defect of the distribution $f_n$, $1-f_n(1)$, represents the probability that an individual at generation $n$ develops this feature. Examples were these models can be applied are discussed in \cite{Iwasa-Michor-Nowak-2004}. For instance, they consider the situation where a virus develops in a vaccinated host and a mutation leads to other type of individuals which are not affected by the vaccine and therefore, they can cause an epidemic outbreak. A second example is a population of cancer cells under some treatment and a mutation might turn them into resistant cells.  The reader is referred to \cite{Iwasa-Michor-Nowak-2004} for more further genetic problems were these models are of interest.

To describe these situations more precisely, let us consider a two-type branching process \linebreak $\{(Y_n(1),Y_n(2))\}_{n\in\N_0}$ defined as follows. The variable $Y_n(1)$ represents the number of individuals of type 1 that are free of the mutation or disease at generation $n$, whereas $Y_n(2)$ denotes the number of individuals of type 2, mutants or individuals with the disease at generation $n$. We assume that individuals of type 1  are the only ones that produce individuals of this type, but they can also give birth to individuals of type 2. Individuals of type 2, however, only produce individuals of the same type. If we also assume independence in the reproduction and that the population starts only with individuals of the first type, then the process is defined as
\begin{align*}
(Y_0(1),Y_0(2))=(N,0),\quad (Y_{n+1}(1),Y_{n+1}(2))=\left(\sum_{i=1}^{Y_n(1)}\xi_{ni}(1),\sum_{i=1}^{Y_n(1)+Y_n(2)}\xi_{ni}(2)\right),\quad n\in\N_0,
\end{align*}
where $N\in\N$, and the variables of the family $\{\xi_{ni}(1),\xi_{ni}(2):i\in\N,n\in\N_0\}$ are independent distributed. If we additionally require  that the distribution of the random vectors $(\xi_{ni}(1), \xi_{ni}(2))$ depend only on $n$,  we are dealing with a (degenerated) two-type branching process in a varying environment. Then, as is easy to see, up to the moment  $\tau_\Delta=\min\{n\in\N:Y_n(2)>0\}$ the first component $\{Y_n(1)\}_{n \in N_0}$ constitutes a DGWPVE, with offspring distributions given by
\begin{align*}
f_n[k]&=P[\xi_{n1}(1)=k,\xi_{n1}(2)=0],\quad k\in\N_0,\\
f_n(1)&=1-\sum_{k=0}^\infty P[\xi_{n1}(1)=k,\xi_{n1}(2)=0]=1-P[\xi_{n1}(2)=0]=P[\xi_{n1}(2)>0].
\end{align*}
For a constant environment this was already observed and applied in \cite{Karlin-Tavare-1982} motivated by a genetic problem in \cite{Robertson-1978}.

Despite the great interest of DGWPVEs and the number of papers dealing with branching processes in a varying environment (see, for instance,  \cite{Jagers-1974}, \cite{Lindvall}, \cite{Agresti}, \cite{Keiding-Nielsen-1975}, \cite{Foster-Goettge-1976}, \cite{Goettge-1976}, \cite{Fujimagari}, \cite{SouzaBigginsa}, \cite{MaCphee-Schuh-1983}, or more recently, \cite{Yu-Pei-2009}, \cite{Hu-Hu-Yin-2011}, and \cite{Kersting-2020}), none of them analyses this type of populations and consequently, our work provides the first results for this problem. Bearing in mind the previous interpretation, we focus on the study of the limiting behaviour of the process. In our first theorem we provide a necessary and sufficient condition for the process to escape the explosion and the absorption (at $0$ or $\Delta$), as occurs for the BPVEs. 
However, we emphasize that the presence of defective distributions in the model together with their change over time (note that we identify time with generation in this framework) makes the difference between their asymptotic behaviour and the aforementioned models. First, we show that unlike the BPVEs, these processes do not become extinct almost surely. Second, we prove that these processes can avoid absorption with a positive probability, contrary to the duality extinction-absorption at $\Delta$ that holds for DGWPs. Indeed, we establish necessary and sufficient conditions for the process to avoid its absorption and we also state the growth rate of the process in the non-absorption set. Our next results describe the asymptotic regimes, as $n\to\infty$, of the probability that the process avoids absorption at generation $n$ and of the expected number of individuals at generation $n$. To determine the size of the process conditionally on its non-absorption at generation $n$, as $n\to\infty$, we introduce the defective branching trees. They are the counterpart of the branching trees in the non-defective case, which have been used in \cite{Geiger-1999} to analyse the limiting behaviour of GWPs, or in \cite{Kerting-Vatutin-2017} for BPVEs. For the readers convenience we explain in full detail how to adapt this construction of branching trees and state their relationship with the DGWPVEs. By this means, we provide an upper bound for the expectation of the process conditionally on its non-absorption at that generation.

%
%



Apart from this introduction, the paper is organised in three sections and one appendix. In Section~\ref{sec:model} we provide the description of the probability model and introduce the notation that we use throughout this paper. In Section~\ref{sec:absorp-explos}, we give several results on the behaviour of the process regarding its absorption either at the state $0$ or $\Delta$. Section~\ref{sec:trees} is devoted to  DGWPVEs conditioned on  non-absorption. In order to ease the reading, we collect the proofs of the results in a final appendix.

\vspace*{0.5cm}

\section{Description of the model}\label{sec:model}

Mathematically, a  {\it defective Galton-Watson process in a  varying environment} $v=\{f_n\}_{n\in\N}$ is a discrete time stochastic process $\{Z_n\}_{n\in\N_0}$ defined recursively as:
\begin{equation}\label{def:model}
Z_0=1,\quad Z_{n}=\sum_{j=1}^{Z_{n-1}}X_{nj},\quad n\in\N,
\end{equation}
where $\{X_{nj}:n, j\in\N\}$ is a family of independent random variables defined on a probability space $(\Omega,\mathcal{A},P)$ and with range $\mathbb N_\Delta=\mathbb N_0\cup\{\Delta\}$,  and for each $n\in\N$, the $X_{nj}$, $j\in\N$, have a common distribution $f_n$ with weights $P[X_{nj}=k]=f_n[k]$, $k\in\N_0$, and generating function
\begin{equation*}\label{eq:fs}
f_n(s)=E\left[s^{X_{n1}}\right]=\sum_{k=0}^\infty f_n[k] s^k,\quad s\in [0,1],
\end{equation*}
satisfying $f_n(1)\leq 1$, thus $P[X_{nj}=\Delta]=1-f_{n}(1)$. Here we use the conventions
\begin{equation}\label{eq:grave-prop}
\Delta+k=\Delta, \quad k\in\N_\Delta, \qquad s^\Delta=0,\quad s\in [0,1],\qquad \sum_{j=1}^\Delta k_j=\Delta, \quad k_j\in\N_\Delta.
\end{equation}
Without further mention in the sequel, we require that \[0< f_n'(1) < \infty\] for all $n \ge 1$.

The possible defect \[\delta_n=1-f_n(1)\] of the distribution $f_n$ is interpreted as the probability with which any particle at generation $n-1$ may send the whole process to an absorbing graveyard state $\Delta$ at generation $n$, where it stays forever. Intuitively, this process is appropriate to represent the evolution of populations where individuals reproduce independently of the others or show certain feature with a positive probability, and the distribution governing the reproduction in each generation is the same for all the individuals. More specifically, the process models populations that evolves as described and that are free of a feature developed by the individuals until the first time that one of them shows it. This event is represented by the state $\Delta$. DGWPVEs have also proved to be a useful tool in analysing other branching processes, as multitype branching processes (see \cite{Braunsteins-Hautphenne-2019}, p.8, where $\Delta$ is identified with the state $\infty$).

The process $\{Z_n\}_{n\in\N_0}$ is an inhomogeneous Markov chain with state space $\N_\Delta$, and two of these states are absorbing: 0 and $\Delta$. Let us denote $f_{k,n}=f_{k+1}\circ\ldots\circ f_n$, for $k=0,\ldots, n$, with the convention that $f_{n,n}(s)=s$, $s\in [0,1]$. From \eqref{def:model}, and by considering the absorption properties \eqref{eq:grave-prop} at the graveyard state, one obtains that $E[s^{Z_n}]=f_{0,n}(s)$, for $s\in [0,1]$.

Let us also denote the probability of extinction and absorption at the graveyard state $\Delta$ by $q=P[Z_n\to 0]$ and $\widehat{q}=P[Z_n\to\Delta]$, respectively, and let $\tau_0=\min\{n\in\N: Z_n=0\}$ be the extinction time, $\tau_\Delta=\min\{n\in\N: Z_n=\Delta\}$ the time of absorption at the graveyard state $\Delta$ and $\tau_a=\tau_0\wedge \tau_\Delta$ the ultimate absorption time. It is straightforward to see that
\begin{equation*}
P[\tau_0\leq n]=f_{0,n}(0),\quad P[\tau_\Delta\leq n]=1-f_{0,n}(1),\quad P[\tau_a\leq n]=1+f_{0,n}(0)-f_{0,n}(1),
\end{equation*}
and
\begin{equation*}
P[\tau_0<\infty]=q,\quad P[\tau_\Delta<\infty]=\widehat{q},\quad P[\tau_a<\infty]=q+\widehat{q}.
\end{equation*}
Observe that if $\widehat{q}=1-q$, then the process becomes absorbed at time $\tau_a$ with probability 1. This was the case for the DGWPs with offspring distribution independent of the generation (see \cite{Sagitov-Minuesa-2017}, p.2); however, as we prove below, the introduction of a varying environment in the model allows the process to avoid absorption and consequently, the explosion of the population might occur. Moreover, note that in the case $f_n(1)=1$, for each $n\in\N$, one obtains the classical branching process in a varying environment. Otherwise, since $\delta_n>0$, for some $n\in\N$, $0<\widehat{q}\leq 1$ and consequently $0\leq q<1$; thus, contrary to the case of Galton-Watson processes in a varying environment (see \cite[Theorem 1]{Kersting-2020}), the almost sure extinction of the population is impossible.

Henceforth, we simply write $E[h(X)]$ to refer $E[h(X); X\neq \Delta]$ for any variable $X\ge 0$ and any function $h(\cdot)$. A main obstacle in the treatment of DGWPVEs compared to the non-defective case is the difficulty to access the moments of $Z_n$. Indeed, we have the formulas
\begin{align}
E[Z_n] &= \prod_{j=1}^n f_j'(f_{j,n}(1)),\label{eq:moments1} \\
\frac{E[Z_n^2]}{E[Z_n]^2}&= \frac 1{E[Z_n]}+ \sum_{j=1}^n \frac{f_j''(f_{j,n}(1))}{f_{j}'(f_{j,n}(1)) \mu_{j,n}},\label{eq:moments2}
\end{align}
with
$$\mu_{j,n}=\prod_{i=1}^j f_i'(f_{i,n}(1)),\quad j=0,\ldots,n,$$
(see \cite[Lemma 4]{Kersting-2020}). These expressions, derived in the Appendix, are in general difficult to control since in the defective case the quantities $f_{i,n}(1)$ can no longer be just  replaced with 1. We come back to this issue.

DGWPVEs can be written in a similar way to controlled branching processes with both the offspring and control distributions depending on the generation as follows:
\begin{equation}\label{eq:model-CBP}
Z_0=1,\quad Z_{n}=\sum_{j=1}^{\phi_{n-1}(Z_{n-1})}\widetilde{X}_{nj},\quad n\in\N,
\end{equation}
where the random variables $\widetilde{X}_{nj}$, $n,j\in\N$, are independent and for each $n\in\N$, the $\widetilde{X}_{nj}$, $j\in\N$, are identically distributed according to
\[g_n(s)=\frac{f_n(s)}{f_n(1)},\quad s\in [0,1].\]
Moreover, the family $\{\phi_n(k):n\in\N_0,k\in\N_\Delta\}$ is a family of independent variables, which is also independent of the family $\{\widetilde{X}_{nj}: n, j\in\N\}$ and such that for each $n\in\N$ and $k\in\N_\Delta$, the random variable $\phi_n(k)$ is defined as
\begin{equation*}
\quad\phi_{n}(k)=\left\{
\begin{array}{llr}
  k&  \text{with probability} & f_{n+1}(1)^k,  \\
 \Delta  & \text{with probability} &1-f_{n+1}(1)^k.
\end{array}
\right.
\end{equation*}
Note that from \eqref{eq:grave-prop}, one has that $\phi_n(\Delta)=\Delta$ a.s., for each $n\in\N_0$. This observation may be used to construct a DGWPVE out of a non-defective process: define the branching process $\{ \widetilde Z_n\}_{n\in\N_0}$ in a varying environment $\{g_n\}_{n\in\N}$,
\[ \widetilde Z_0=1 , \qquad \widetilde Z_n = \sum_{j=1}^{\widetilde Z_{n-1}} \widetilde X_{nj}, \qquad n \in \N,\]
with values in $\N_0$, and the stopping time
\[ \tau_\Delta = \min \{ n\in \mathbb N: \phi_{n-1}(\widetilde Z_{n-1})= \Delta\} .\]
Then
\begin{align}\label{X}
Z_n = \begin{cases} \tilde Z_n , &\text{if }  \tau_\Delta>n, \\ \Delta, &\text{if } \tau_\Delta \le n,\end{cases}
\end{align}
is a DGWPVE $v=\{f_n\}_{n\in\N}$. This coupling will be useful in different proofs below.

%

\vspace*{0.5cm}

\section{Some results on absorption and explosion}\label{sec:absorp-explos}

In this section, we provide some results related to the absorption and explosion of a DGWPVE. The first result establishes the almost sure convergence of the process $\{Z_n\}_{n\in\N_0}$ to a random variable $Z_\infty$ with values in the set $\N_\Delta \cup \{\infty\}$. As mentioned before, for a DGWPVE the duality explosion-absorption (either at the state 0 or $\Delta$) does not always hold true, and similarly to the Galton-Watson process in varying environment (see Lindvall \cite{Lindvall}), the event $\{0<Z_\infty<\infty\}$ may have positive probability; the next theorem also establishes  for that situation the same necessary and sufficient condition as was detected by  Lindvall in the non-degenerate case.

\begin{thm}\label{thm:Z-infty}
Let $\{Z_n\}_{n\in\N_0}$ be a DGWPVE $v=\{f_n\}_{n\in\N}$. There exists a random variable $Z_\infty$ with values in $\N_\Delta \cup \{\infty\}$ such that as $n\to\infty$,
\begin{equation*}
Z_n\to Z_\infty\quad \text{a.s.}
\end{equation*}
Moreover,
\begin{equation}\label{eq:thm-Z-infty}
P[0<Z_\infty<\infty]>0\qquad\Leftrightarrow\qquad \sum_{n=1}^\infty (1-f_n[1])<\infty.
\end{equation}
\end{thm}

Observe that from Theorem \ref{thm:Z-infty}, one has that a necessary condition for the event $\{0<Z_\infty<\infty\}$ to have positive probability is that  $f_n[1]\to 1$, as $n\to\infty$ (and consequently, $\delta_n\to 0$); indeed, those convergences need to be fast in the terms described in \eqref{eq:thm-Z-infty}. These facts are illustrated in the examples below. Let us remark that the necessary and sufficient condition provided in Theorem \ref{thm:Z-infty} is equivalent to the notion of complete convergence of the offspring random variables $X_{n1}$ to 1, as introduced by Hsu and Robbins \cite{Hsu-Robbins-1947} and meaning that $\sum_{n=1}^\infty P[|X_{n1}-1|>\epsilon]<\infty$ for all $\epsilon >0$. Indeed, if $\epsilon\in (0,1)$, then 
$$\sum_{n=1}^\infty (1-f_n[1])=\sum_{n=1}^\infty P[|X_{n1}-1|>\epsilon].$$
Taking into account that the variables of the family $\{X_{n1}: n\in\N\}$ are independent, by using the first and second Borel Cantelli lemmas, we  have that the previous condition is equivalent to the almost sure convergence $X_{n1}\to 1$ as $n\to\infty$. This renders some plausibility to  the equivalence \eqref{eq:thm-Z-infty}. Let us mention that the somewhat stronger condition $\sum_{n=1}^\infty E[|X_{n1}-1|]<\infty$ has been used by Isp\'any \cite{Ispany-2016} in the context of branching processes {\em with immigration} in a varying environment.

\begin{ex}
\begin{enumerate}[label=(\alph*),ref=\emph{(\alph*)}]
\item Let us consider a DGWPVE with environment satisfying that $f_1(s)=2^{-1}s$ and $f_n(s)=(1-n^{-1})s$, for $n\geq 2$ and $s\in [0,1]$. Then,
$$\sum_{n=1}^\infty (1-f_n[1])=\sum_{n=1}^\infty \delta_n=\frac{1}{2}+\sum_{n=2}^\infty\frac{1}{n}=\infty,$$
and consequently, $P[0<Z_\infty<\infty]=0$.
\item Let us consider a DGWPVE such that for $s\in [0,1]$, $f_1(s)=2^{-1}s$ and $f_n(s)=(1-n^{-2})s$, for $n\geq 2$. Then,
$$\sum_{n=1}^\infty (1-f_n[1])=\sum_{n=1}^\infty \delta_n=\frac{1}{2}+\sum_{n=2}^\infty\frac{1}{n^2}<\infty,$$
and consequently, $P[0<Z_\infty<\infty]=P[Z_\infty=1]>0$.
\end{enumerate}
\end{ex}

\vspace*{0.5cm}

In view of Theorem \ref{thm:Z-infty}, it is interesting to study the conditions for having the duality extinction-absorption at $\Delta$ that holds for DGWPs; Theorem~\ref{thm:explosion} and Proposition~\ref{thm:explosion-2} deal with that issue. Let us first introduce the notation:
$$ \mu_0=1,\qquad \mu_n=\prod_{i=1}^n f_i'(1),\quad n\in\N.$$

\vspace{2ex}


\begin{thm}\label{thm:explosion}
Let $\{Z_n\}_{n\in\N_0}$ be a DGWPVE $v=\{f_n\}_{n\in\N}$. Assume that there exists $c>0$ such that for all $n\in\N$,
\begin{equation}\label{eq:cond}
E[X_{n1}^2;X_{n1}\geq 2]\leq c E[X_{n1};X_{n1}\geq 2]E[X_{n1}|X_{n1}\geq 1].
\end{equation}
Then
\begin{equation*}
P[\tau_a=\infty]>0\qquad\Leftrightarrow\qquad \inf_{n\in\N} \mu_n >0,  \quad\sum_{n=1}^\infty \delta_n \mu_{n-1} < \infty \quad \text{and}\quad \sum_{n=1}^\infty \frac{f_n''(1)}{f_n'(1)\mu_{n}}<\infty .
\end{equation*}
\end{thm}


Obviously the regularity condition in \eqref{eq:cond} is  equivalent to
\begin{equation}\label{eq:cond_2}
E[\widetilde X_{n1}^2;\widetilde X_{n1}\geq 2]\leq c E[\widetilde X_{n1};\widetilde X_{n1}\geq 2]E[\widetilde X_{n1}|\widetilde X_{n1}\geq 1].
\end{equation}
It is a uniformity assumption introduced in \cite{Kersting-2020} for the study of branching processes in a varying environment and it is satisfied by a wide class of  distributions.

%
%

\begin{ex} As an application of the previous theorem, in the following examples we show two  DGWPVEs in which the explosion set has null probability and positive probability, respectively.

\begin{enumerate}[label=(\alph*),ref=\emph{(\alph*)}]
\item Let us consider a DGWPVE with environment satisfying that for $s\in [0,1]$, $f_n(s)=(1-n^{-1}2^{-n})s^2$, $n\in \N$. It is easy to check that condition \eqref{eq:cond} holds for $c\geq 2$ and since $f_n[1]=0$ for $n\in\N$, $P[0<Z_\infty<\infty]=0$. Moreover, due to the fact that $f_n[0]=0$, $n\in\N$, one has $P[Z_n\to 0]=0$, and since $\prod_{i=1}^n (1-i^{-1}2^{-i})$ converges to a non-zero number,
$$\sum_{n=1}^\infty \delta_n\mu_{n-1}=\frac{1}{2}\sum_{n=1}^\infty \frac{1}{n}\prod_{i=1}^{n-1} \bigg(1-\frac{1}{i2^{i}}\bigg)=\infty,$$
implying $P[\tau_a=\infty]=0$. In particular, this process gets absorbed at $\Delta$ a.s.

\item Let us consider a DGWPVE such that for $s\in [0,1]$, $f_n(s)=(1-n^{-2}2^{-n})s^2$, $n\in \N$. With the same arguments as in the previous example, it is easy to verify that $P[0<Z_\infty<\infty]=0$ and $P[Z_n\to 0]=0$, and condition \eqref{eq:cond} holds true for $c\geq 1$. Now, we use the fact that $\prod_{i=1}^n (1-i^{-2}2^{-i})$ converges to a non-zero number to get that $\mu_n\to\infty$, as $n\to\infty$,
$$\sum_{n=1}^\infty \delta_n \mu_{n-1}=\frac{1}{2}\sum_{n=1}^\infty \frac{1}{n^2}\prod_{i=1}^{n-1} \bigg(1-\frac{1}{i^22^{i}}\bigg)<\infty,$$
and
$$\sum_{n=1}^\infty \frac{f_n''(1)}{f_n'(1)\mu_n}=\sum_{n=1}^\infty \frac{1}{2^n}\prod_{i=1}^n\frac{1}{1-\frac{1}{i^22^i}}<\infty,$$
and consequently, $P[\tau_a=\infty]=P[Z_\infty=\infty]>0$.
\end{enumerate}
\end{ex}

\vspace{0.5cm}

The coupling in \eqref{X} also enables to describe the asymptotic behaviour of the defective process $\{Z_n\}_{n\in\N_0}$ on the non-absorption set $\{\tau_a=\infty\}$. Observe that the sequence of random variables $\{\widetilde{W}_n\}_{n\in\N_0}$, with $\widetilde{W}_n=\widetilde{Z}_n\prod_{i=1}^n \frac{f_i(1)}{f_i'(1)}$, is a non-negative martingale, and hence, there exists an integrable and non-negative random variable $\widetilde{W}$ such that $\widetilde{W}_n$ converges a.s. to $\widetilde{W}$, as $n\to\infty$.  Moreover, under the equivalent conditions formulated in  Theorem \ref{thm:explosion} it follows that $\sum_{n = 1}^\infty\delta_n< \infty$, thus $\prod_{i=1}^n f_i(1)= \prod_{i=1}^n(1-\delta_i)$ converges to a strictly positive limit and then, $\widetilde{Z}_n/\mu_n$ converges to a non-negative random variable $W$. Also, because of \eqref{eq:cond_2} we have from \cite[Theorem 2]{Kersting-2020} that $\{ \widetilde Z_\infty>0\}=\{\widetilde{W}>0\} = \{W>0\}$ a.s. The following result is straightforward from the fact that on the set $\{\tau_a=\infty\}$, both processes $\{Z_n\}_{n\in\N_0}$ and $\{\widetilde{Z}_n\}_{n\in\N_0}$ visit the same states and consequently $\{\tau_a=\infty\}\subset\{\widetilde{Z}_\infty>0\}=\{W>0\}$ a.s.

\begin{cor}\label{cor:asympt-Z-infty}
Under the assumption \eqref{eq:cond} from Theorem \ref{thm:explosion}, if $P[\tau_a=\infty]>0$, then there exists a non-negative random variable $W$ such that, as $n\to\infty$,
\begin{equation*}
\frac{Z_n}{\mu_n}\to W>0\quad\text{ a.s. on } \{\tau_a=\infty\}.
\end{equation*}
\end{cor}


Theorem \ref{thm:explosion} shows that at least $\delta_n \to 0$ as $n \to \infty$ is required in order to have a strictly positive  probability $P[\tau_a=\infty]$ of non-absorption. Conversely, in the situation, where e.g. $\inf_{n \ge 1}\delta_n>0$ holds, the question arises how to estimate the absorption probabilities $P[\tau_a>n]$.
Applying the mean-value theorem to the formula $P[\tau_a>n]= f_{0,n}(1)-f_{0,n}(0)$ yields
\[ \prod_{i=1}^n f_i'(f_{i,n}(0)) \le P[\tau_a>n] \le \prod_{i=1}^n f_i'(f_{i,n}(1)),\quad n\in\N . \]
The lower bound can be 0 and the upper bound greater than or equal to 1, therefore one may ask whether one can do better. In a first step, we provide the following result.

\begin{prop}\label{thm:explosion-2}
For each $n\in\N$,
\[ P[\tau_a >n] \le \inf_{i\le n} \mu_i.\]
Moreover,  if there exists a constant $c>0$ such that for all $n\in\N$,
\begin{equation}\label{cond-Xni}
E[X_{n1}^2; X_{n1}\ge 2] \le c E[X_{n1}; X_{n1} \ge 2] ,
\end{equation}
then, there is a constant $c'>0$ such that as $n \to \infty$
\[ \frac{E[Z_n]^2}{E[Z_n^2]} \le P[\tau_a >n] \le c'\frac{E[Z_n]^2}{E[Z_n^2]} .\]
\end{prop}

Condition \eqref{cond-Xni} is  stronger than \eqref{eq:cond_2}, it implies that $E[X_{n1}^2]$ is bounded uniformly in $n$ (see formula \eqref{fnprime} below).

Combining these bounds on $P[\tau_a >n]$  with the equation \eqref{eq:moments2} confronts us with the problem to obtain  estimates for the terms  $f_{i,n}(1)$ from the given distributions $f_n$. In general this appears to be difficult. The following theorem considers the case where all solutions $\theta_n$   of the equations $f_n(s)=s$  belong to some interval $[ \rho, \sigma]$ with $0<  \rho< \sigma < 1$. For $0< s\le1$ denote
\[ \mu_n(s) = \prod_{i=1}^n f_i'(s), \quad \nu_n(s)= \sum_{i=1}^n \frac 1{\mu_i(s)}.\]

\begin{thm}\label{cor:asympt-E-infty}
Under the assumption \eqref{cond-Xni}, then for any $ \rho \in (0,1)$ and $n_0\in \mathbb N$ such that $f_n( \rho)\ge  \rho$ for all $n \ge n_0$ we have
\[ \liminf_{n\to \infty} \frac{ E[Z_n]}{\mu_n( \rho)} >0, \quad \text{and}\quad \liminf_{n \to \infty} \nu_n( \rho) P[\tau_a >n] >0  .\]
If, in addition, $\inf_{n \ge 1}\delta_n>0$, $\inf_{n \ge 1} f_n''( \rho)>0$,  then for any $ \sigma\in ( \rho,1)$  such that $f_n( \sigma)\le  \sigma$ for all $n \ge n_0$, and for any $\varepsilon >0$ we have
\[ \limsup_{n\to \infty} \frac{ E[Z_n] }{\mu_n( \sigma+ \varepsilon)}< \infty, \quad \text{and}\quad \limsup_{n \to \infty} \nu_n( \sigma+ \varepsilon) P[\tau_a >n] <\infty .\]
\end{thm}

\medskip

\begin{ex}
\begin{enumerate}[label=(\alph*),ref=\emph{(\alph*)}]
\item If $\inf_{n \ge 1} f_n[0] >0$ then we may set $\rho=\inf_{n \ge 1} f_n[0]$. If $\sup_{n \ge 1} f_n[1] <1$ then we may set $\sigma= \sup_{n \ge 1} f_n[1]$.
\item   In some important cases the numbers $\rho$ and $\sigma$ can be explicitly evaluated. First let us look at the case of defective {\bf binary offspring} meaning that $f_n[k]=0$ for $k >2$ and all $n \ge 1$. Denote $p_n= f_n[2]$, $q_n=f_n[1]$ and $r_n=f_n[0]$ and suppose that there are $0<\theta_n<1$ fulfilling
\[ f_n(\theta_n)= \theta_n.\]
These equations are quadratic with two roots,  and the $\theta_n$  is the smaller root, given by
\[ \theta_n = \frac {1-q_n}{2p_n}- \sqrt{\frac{(1-q_n)^2}{4p_n^2 }-
\frac{r_n}{p_n}} .\]
It follows
\[ \rho= \inf_{n \ge 1} \theta_n  , \quad \sigma= \sup_{n \ge 1} \theta_n,\]
 whenever $\inf_{n \ge 1}\delta_n>0$ and $\inf_{n \ge 1}r_n>0$.

Second we consider the case of defective {\bf linear fractional} distributions. In this case the generating functions are of the form
\[ f_n(s)= q_n + \frac{r_n}{1-p_ns} , \]
and the equation $f_n(\theta_n)=\theta_n$ with $0<\theta_n < 1$ reduces once more to a quadratic equation resulting in
\[ \theta_n = \frac{1+p_nq_n}{2p_n} - \sqrt{ \frac{(1+p_nq_n)^2}{4p_n^2}- \frac{r_n+q_n}{p_n}}. \]
Again, $\rho$ and $\sigma$ are the infimum and supremum of these quantities, whenever $\inf_{n \ge 1}\delta_n>0$ and $\inf_{n \ge 1}(q_n+r_n)>0$.

\item Let $f$ be a distribution with defect $1-f(1)>0$ and let $f_n$ be a sequence of defective distributions such that $f_n(s) \to f(s)$ as $n \to \infty$ for all $0\le s \le 1$. Let $\theta$ be the unique solution of the equation  $f(\theta)=\theta$ and assume that $\theta >0$. Then, letting $ \rho=\theta-\eta$ and  $ \sigma= \theta+\eta$, all assumptions of Theorem \ref{cor:asympt-E-infty} are satisfied for any  $0<\eta<\min\{\theta,1-\theta\}$.

Also, since $f'_n(s)\to f'(s)$ for all $s<1$ and since $f'(s)$ is strictly increasing, it follows that for any $\eta >0$
\[  \lim_{n \to\infty} \frac{\mu_n(s)}{f'(s-\eta)^n} =\infty  \quad \text{and} \quad \lim_{n \to \infty} \frac{\mu_n(s)}{f'(s+\eta)^n} =0 ,\]
or equivalently
\[ \lim_{n \to \infty} \frac 1n \log \mu_n(s)= \log f'(s). \]

%

If $f'(s)<1$, we observe that $\mu_n(s)$ is exponentially decreasing in $n$ at the negative rate $\log f'(s)$,  and then, as is easy to see,  $\nu_n(s)$ is increasing at the rate $-\log f'(s)$, that is
\[ \lim_{n \to \infty} \frac 1n \log \nu_n(s)= - \log f'(s). \]

Now $f'(\theta)<1$. Therefore, applying Theorem \ref{cor:asympt-E-infty} we obtain
\[ \lim_{n \to \infty} \frac 1n \log E[Z_n] = \log f'(\theta) \quad \text{and} \quad  \lim_{n \to \infty} \frac 1n \log P[\tau_a>n] = - \log f'(\theta) .\]

\end{enumerate}
\end{ex}


\vspace*{0.5cm}

\begin{rem}
Along similar lines we can get estimates on  extinction probabilities. Let
\[ q_l = \lim_{n \to \infty} f_{l,n}(0) \]
which is the probability of ultimate extinction given that there is one individual at generation $l$. Note that $q_{l}=f_{l,n}(q_n)$ for $0<l\le n$. Let us suppose again that there is a number $\sigma \in (0,1)$ such that
\[f_n(\sigma) \le \sigma\] for all $n \ge 1$. Then, along the same lines as above, it follows that $f_{l,n}(0) \le f_{l,n}(\sigma)\le \sigma$ and consequently
\[ q_l \le \sigma\]
for all $l\ge 1$. Now
\begin{align*}
P[n<\tau_0 <\infty] = q_0- f_{0,n}(0) = f_{0,n}(q_n)- f_{0,n}(0) \le f_{0,n}(\sigma) -f_{0,n}(0) .
\end{align*}
Applying the mean-value theorem we have $f_{0,n}(\sigma) -f_{0,n}(0) \le \sigma \prod_{i=1}^n f_{i}'(f_{i,n}(\sigma))$ and finally
\[  P[n<\tau_0 <\infty] \le \sigma\prod_{i=1}^n f_{i}'(\sigma). \]
Similarly one derives the lower bound $ P[n<\tau_{\Delta}<\infty]\ge (1-\sigma)\prod_{i=1}^n f_{i}'(f_{i,n}(\sigma))$.
\end{rem}

\vspace*{0.5cm}

\section{Conditioned DGWPVEs and defective branching trees}\label{sec:trees}

In this section we consider the random variables $Z_n$, conditioned on the event $\{\tau_a >n\}$. In particular, we obtain the following theorem on the conditional expectation of $Z_n$, which can be considered as a counterpart of Theorem 2 (a) in \cite{Sagitov-Minuesa-2017} in the case of varying environments. For other conditional moments similar estimates are also feasible but we omit them for sake of brevity.

\begin{thm}\label{thm:cond-expectation}
If $\inf_{n\in\N}f_n(0)=\alpha>0$, and $\sup_{n\in\N}f_n(1)=\beta<1$, then for some $c>0$
\begin{align*}
E[Z_n \mid \tau_a>n]\leq   1+cf_n'(1)\sum_{j=0}^{n-1}\beta^{j}\left(1+\frac{f_{n-j}''(1)}{f_{n-j}'(1)}\right).
\end{align*}
In particular, if condition \eqref{cond-Xni} also holds, then the sequence $\{E[Z_n \mid \tau_a>n]\}_{n\in\N_0}$ is bounded.
\end{thm}

It seems  little promising to prove this result by those tools we used so far, at least Theorem \ref{cor:asympt-E-infty} provides for this purpose  insufficient control on the probability of the event $\{\tau_a>n\}$. Therefore we come up with a different means, namely a probabilistic construction of the
 {\it conditional defective branching tree.}  Such branching trees have been applied to study properties of GWPs (see \cite{Geiger-1999}) or BPVEs (see \cite[Chapter 2]{Kerting-Vatutin-2017}). These structures are richer frameworks than the family of variables $\{X_{ni}:n,i\in\N\}$ associated with the branching process since they provide information on the genealogical relationship between any two individuals of the population.

\vspace*{0.25cm}

\subsubsection*{Defective family trees}

For the readers convenience we introduce this approach in full detail. First, we introduce the formal definitions of  {\it defective family tree}. To that end, we make use of the Ulam-Harris labelling with some modifications in order to indicate the absorption at the state $\Delta$. Let us consider the two following sets
\begin{equation*}
\I=\bigcup_{n=0}^\infty \N^n,\quad\text{ and }\quad \dI=\bigcup_{n=0}^\infty (\N\cup\{\Delta\})^n,
\end{equation*}
with the convention $\N^0=\varnothing$ and $(\N\cup\{\Delta\})^0=\varnothing$. Individuals in the population correspond to elements of the set $\I$ written as finite strings, where the element $\varnothing$ is known as the {\it founding ancestor} or the {\it root}. Moreover, given an individual $i\in\I$, we denote $g(i)$ the generation of the individual $i$ and write $i\Delta$ to indicate that the individual $i$ sends the process to the graveyard state $\Delta$ at the next generation. Thus, the elements of type $i\Delta$ do not represent individuals of the population and we refer to them as {\it defective elements}. We also write $c(i)$ to refer to the number of children of the individual $i$ when $i$ does not send the process to $\Delta$ (including the case $c(i)=0$), and write $c(i)=\Delta$ otherwise. Moreover, we make the following assumptions:
\begin{itemize}
\item If $g(i)=n$, then $i$ is a string $j_1\ldots j_n$, with $j_1,\ldots,j_n\in\N$.
\item If $i$ has $c=c(i)\in\N$ children, then the label of each child is $ij$, $1\leq j\leq c$.
\item If $n\in\N$, then the progenitor of $i=j_1\ldots j_n$ is the truncated string $i'=j_1\ldots j_{n-1}$.
\end{itemize}

\vspace*{0.25cm}

Under the previous considerations, the population is included in a subset $t$ of $\dI$, whose elements satisfy the following properties:
\begin{enumerate}[label=(\roman*),ref=(\roman*)]
\item $\varnothing \in t$.\label{cond:founding-ancestor}
\item Let $i \in \I$, and $j \in \N_\Delta$. If $i j \in t$, then $i \in t$. 
\item If there exists $i \in t$ equal to a string $j_1\ldots j_n\Delta$, with $j_1,\ldots, j_n\in\N$, then all the elements of $t$ are strings of length at most $n+1$. 
\item Let $i \in \mathcal{I}$ and $j \in \mathbb{N}_\Delta$. If $i j \in t$, then either $j\in\N$ and $i j' \in t$ for all $1 \leq j' \leq j$, or else $j=\Delta$ and $ij'\not\in t$ for all $j'\in\N$. 
\item For $i \in t$, there exists a $j \in \mathbb{N}$ such that $i j \notin t$. \label{cond:finite-children}
\end{enumerate}
A set $t$ satisfying conditions \ref{cond:founding-ancestor}-\ref{cond:finite-children} is called {\it defective family tree}.


\begin{figure}[H]
\centering\includegraphics[width=0.4\linewidth]{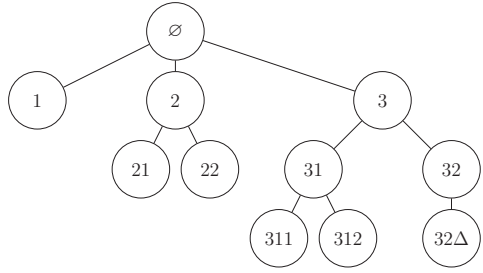}
\caption{Example of a defective family tree with a defective element at the third generation.  The height is $h(t)=2$ and the population sizes are $Z_1(t)=3$, $Z_2(t)=4$, and $Z_3(t)=\Delta$.}\label{fig:def-family-tree}
\end{figure}

We now define the height and generation sizes associated with a defective family tree. To that end, let us introduce the following sets
\begin{equation*}
\I^{(n)}=\bigcup_{l=0}^n \N^l,\quad\text{ and }\quad \dI^{(n)}=\bigcup_{l=0}^n (\N\cup\{\Delta\})^l,
\end{equation*}
Then, the {\it height of the defective family tree $t$} is
\begin{align*}
h(t)&=
\max\{n\ge 0:t\cap\dI^{(n)}\subseteq \I^{(n)} \text{ and } t\cap\I^{(n)} \neq t\cap\I^{(n-1)} \},
\end{align*}
and it can be finite or infinite. That is: if the population gets extinct or is absorbed at generation $n$, then the height of the tree is $n-1$, and else the height is infinite. The {\it generation sizes of the tree $t$} are defined as follows:
\begin{align*}
z_n(t)&=\left\lbrace
\begin{array}{ll}
\sharp\{i\in t:g(i)=n\}, &\quad \text{if } n\leq h(t),\\
\Delta, &\quad \text{if } n> h(t)\text{ and there exists }j_1\ldots j_{h(t)}\Delta\in t,\\
0, &\quad \text{if } n> h(t)\text{ and }j_1\ldots j_{h(t)}\Delta\notin t,\text{ for any }j_1,\ldots, j_{h(t)}\in\N,
\end{array}
\right.
\end{align*}
(see Figure~\ref{fig:def-family-tree} for illustration).


Let $t,t'$ be two defective family trees and $h\in\N_0$, we write $t\stackrel{h}{=}t'$ if and only if
$$t\cap \dI^{(h)}=t'\cap \dI^{(h)},$$
and we say that $t$ and $t'$ coincide up to the height $h$. 
Let us denote the set of all the defective family trees $t\subseteq \dI$ by $\T$ and consider the $\sigma$-field generated by all the sets of the form
$$\{t'\in\T:t\stackrel{h}{=}t'\},\quad h\in\N_0,\ t\in\T$$
Then,
$t\mapsto h(t)$ and $t\mapsto z_n(t)$ are measurable mappings for $n\in\N_0$. 

\vspace*{0.25cm}

\subsubsection*{Defective branching trees in varying environment}

Now, we introduce the definition of defective branching trees in a varying environment emulating those for the non-defective case in \cite[Chapter 2]{Kerting-Vatutin-2017}. A {\it defective branching tree in varying environment (DBTVE) $v=\{f_n\}_{n\in\N}$} is a $\T$-valued random variable $T$ with distribution given as follows:
\begin{enumerate}
\item[(i)] If $h\in\N_0$ and $t\subseteq \I$, then
$$P_v\left[T\stackrel{h}{=}t\right]=\prod_{i\in t\cap\I:g(i)<h} f_{g(i)+1}[c(i)].$$
\item[(ii)] If $h\in\N_0$, $t\in\T$ is such that $t\nsubseteq\I$ and $h\leq h(t)$, then
$$P_v\left[T\stackrel{h}{=}t\right]=\prod_{i\in t\cap\I:g(i)<h} f_{g(i)+1}[c(i)].$$
\item[(iii)] If $h\in\N_0$, $t\in\T$ is such that $t\nsubseteq\I$ and $h>h(t)$, then
\begin{align*}
P_v\left[T\stackrel{h}{=}t\right]&=\prod_{i\in t\cap\I:g(i)<h(t)-1} f_{g(i)+1}[c(i)]\prod_{\substack{i\in t\cap\I:g(i)=h(t),\\ c(i)\neq\Delta}} f_{g(i)+1}[c(i)]\cdot\\
&\phantom{=}\cdot\prod_{\substack{i\in t\cap\I:g(i)=h(t),\\ c(i)=\Delta}} (1-f_{g(i)+1}(1)).
\end{align*}
\end{enumerate}

\vspace*{0.5cm}

Given a DBTVE $T$, we can define a DGWPVE $v=\{f_n\}_{n\in\N}$ associated with $T$; this is the process $\{Z_n\}_{n\in\N_0}$ given by
\begin{align*}
Z_n(\omega)= z_n(T(\omega)), \quad \text{ for each }\omega\in\Omega,n\in\N_0.
\end{align*}



By using the relation between the DBTVE and the DGWPVE associated with it, we define the extinction time, absorption time at $\Delta$, and ultimate absorption time of a DBTVE as the corresponding absorption times in the associated DGWPVE. Again, these are denoted as
$\tau_0=\min\{n\in\N:Z_n=0\}$, $\tau_\Delta=\min\{n\in\N:Z_n=\Delta\}$, $\tau_a=\min\{\tau_0,\tau_\Delta\}$.
 Then, it is immediate that for $j,n\in\N_0$,
\begin{align*}
P_{v_j}[\tau_0\leq n]&=f_{j,n+j}(0),\\
P_{v_j}[\tau_\Delta\leq n]&=1-f_{j,n+j}(1),\\
P_{v_j}[\tau_a\leq n]&=1-f_{j,n+j}(1)+f_{j,n+j}(0),
\end{align*}
where $v_j=\{f_{j+1},f_{j+2},\ldots\}$ (note that $v=v_0$). We write $\tau_0(T)$, $\tau_\Delta(T)$ and $\tau_a(T)$, if we like to point out that the absorption times stem from the random tree $T$.


Let us also denote the subtrees of the random tree $T$ founded by the children of the root as  $T_i$, $1\leq i\leq z_1(T)$, and let $R_{n}$ be the rank of the left-most child of the root that has a descendant in generation $n$,
$$R_{n}=\min\left\{1\leq i\leq z_1(T):0<z_{n-1}(T_i)\neq \Delta\right\},\quad n\in\N_0, $$
with $\min\varnothing=\infty$. We note that $R_{n}<\infty$ a.s. on the set $\{\tau_a>n\}$.
Similarly to DBTVEs, a defective tree $t$ satisfying $h(t)\geq 1$ can be divided into $c$ subtrees $t_1,\ldots,t_c$, where $c=c(\varnothing)$ is the number of offspring of the founding ancestor. Moreover, these subtrees can be written as $t_j=\{i\in\dI:ji\in t\}$, for $j=1,\ldots,c$.

Now we like to investigate defective branching trees conditioned to achieve a certain height. They can be built up by a stepwise probabilistic construction as introduced by Geiger \cite{Geiger-1999}  for standard Galton-Watson processes. The single steps proceed according to the following lemma, which gives insight into the structure of such conditional defective trees.

\begin{lem}\label{lem:construction-branching-tree}
Let $n\in\N$ and $v=\{f_n\}_{n\in\N}$ be a varying environment. Let the random defective tree $T$ be composed of the defective subtrees $T_1,\ldots,T_C$ in generation one, where $C$ is the random variable denoting the number of children of the founding ancestor and let $D$ be a $\N$-valued random variable with $1\leq D\leq C$. Assume that:
\begin{enumerate}
\item[\emph{(i)}] The distribution of $(D,C)$ is a proper probability distribution given by $\{g_n[d,c]\}_{1\leq d\leq c;c\in\N}$, with
$$g_n[d,c]=\frac{f_{1,n}(1)-f_{1,n}(0)}{f_{0,n}(1)-f_{0,n}(0)}f_1[c]f_{1,n}(0)^{d-1}f_{1,n}(1)^{c-d}.$$
\item[\emph{(ii)}] Given $(D,C)$, the random trees $T_1,\ldots,T_C$ are independent and satisfy that
\begin{enumerate}[label=\emph{(\alph*)},ref=\emph{(\alph*)}]
\item[\emph{(a)}] $T_i$, $1\le i < D$, are DBTVEs $v_1=\{f_2,f_3,\ldots\}$ conditioned on $\{\tau_0(T_i)\le n-1\}$.
\item[\emph{(b)}] $T_D$ is a DBTVE $v_1=\{f_2,f_3,\ldots\}$ conditioned on $\{\tau_a(T_D)>  n-1\}$.
\item[\emph{(c)}] $T_i$,  $D<i\le C$, are DBTVEs $v_1=\{f_2,f_3,\ldots\}$ conditioned on $\{\tau_\Delta(T_i)> n-1\}$.
\end{enumerate}
\end{enumerate}
Then, $T$ is a DBTVE $v$ conditioned on $\{\tau_a>n\}$, that is, to have height at least $n$. 
\end{lem}

%

Given a defective tree $t$ with height $h(t)\geq l$, then we call {\it distinguished path} or {\it spine of length $l$} to a sequence $d_1,\ldots,d_l\in\N$ such that any string of length $k\leq l$ satisfies $d_1\ldots d_k\in t$. Analogously, a {\it distinguished path of length $n$} in a defective branching tree $T$ results from the numbers $D_1,\ldots,D_n$ and the distinguished individuals $\Lambda_1,\ldots,\Lambda_n$, with $\Lambda_l=D_1\ldots D_l$, $l=1,\ldots,n$. It starts with the founding ancestor and it finishes in the distinguished individual $\Lambda_n$ at generation $n$.

\vspace{2ex}

The previous lemma allows us to construct a DBTVE conditioned on non-absorption up to height $n$ as indicated below.

\vspace{2ex}

\begin{construction}\label{construc:tree}
Let us consider the probability distributions $\{g_{l,n}[d,c]\}_{1\leq d\leq c;c\in\N}$, with $l=1,\ldots,n$ and
$$g_{l,n}[d,c]=\frac{f_{l,n}(1)-f_{l,n}(0)}{f_{l-1,n}(1)-f_{l-1,n}(0)}f_l[c]f_{l,n}(0)^{d-1}f_{l,n}(1)^{c-d},\quad d=1,\ldots,c; c\in\N.$$
Then, the following procedure provides a DBTVE $v=\{f_n\}_{n\in\N}$ conditioned upon the event $\{\tau_a>n\}$:
\begin{enumerate}
\item Establish the distinguished path containing the founding ancestor $\varnothing$ and $n$ distinguished individuals $\Lambda_1,\ldots,\Lambda_n$ from generation 1 to $n$.

\item Generate independent random variables $(D_1,C_1),\ldots,(D_n,C_n)$ with distributions $g_{1,n},\ldots,\linebreak g_{n,n}$. Supply the distinguished individual in generation $l$, $\Lambda_l$, with $C_l-1$ siblings, $D_l-1$ of them to the left and $C_l-D_l$ to the right.

\item Given these random variables, generate independent defective branching trees $T_{1,l},\ldots,\linebreak T_{D_l-1,l},T_{D_l+1,l},\ldots,T_{C_l,l}$ in varying environment $v_l=\{f_{l+1},\ldots\}$ satisfying:
\begin{enumerate}
\item $T_{i,l}$, $1 \le i <D_l$, are DBTVEs conditioned on $\{\tau_0(T_{i,l})< n-l\}$,
\item $T_{i,l}$, $D_l<i\le C_l$, are DBTVEs conditioned on $\{\tau_\Delta(T_{i,l})\ge n-l\}$.
\end{enumerate}
Attach the trees $T_{1,l},\ldots,T_{D_l-1,l}$ to the siblings to the left of the distinguished path and  $T_{D_l+1,l},\ldots,T_{C_l,l}$ to the siblings to the right in generation $l$.
\item Complete the tree by adding an independent, unconditioned defective branching tree $T_{D_n,n}$ in varying environment $v_n=\{f_{n+1},\ldots\}$ on the top of the distinguished individual in generation $n$.
\end{enumerate}
\end{construction}

\vspace*{0.25cm}

Then, by this previous construction, an induction argument and Lemma \ref{lem:construction-branching-tree} lead to the following result.
\begin{prop}\label{prop:construct-tree-contion}
Let $\widetilde{T}_n$ be the branching tree resulting from the previous construction. Then:
$$\mathcal{L}(\widetilde{T}_n)=\mathcal{L}(T|\tau_a(T)>n),\quad n\in\N_0,$$
where $T$ is a DBTVE $v=\{f_n\}_{n\in\N}$.
\end{prop}

\begin{appendix}
%
%

%
\section*{Appendix}

Before providing the proofs of the results of the paper, we derive the expressions \eqref{eq:moments1} and \eqref{eq:moments2}. The first one follows immediately from the fact that $E[Z_n]=f_{0,n}'(1)$ and by recursion we get
$$E[Z_n]=f_{0,n}'(1)=f_{0,n-1}'(f_n(1))f_n'(1)=\ldots=\prod_{j=1}^n f_{j}'(f_{j,n}(1))=\mu_{n,n}.$$
Next, to show the expression in \eqref{eq:moments2} we write 
\begin{align*}
\frac{E[Z_n^2]}{E[Z_n]^2}&=\frac{1}{E[Z_n]}+\frac{E[Z_n(Z_n-1)]}{E[Z_n]^2}= \frac{1}{E[Z_n]}+ \frac{f_{0,n}''(1)}{f_{0,n}'(1)^2},
\end{align*}
and use the fact that
\begin{align*}
\frac{f_{0,n}''(1)}{f_{0,n}'(1)}&= [\log(f_{0,n}'(1))]' = \left[\sum_{j=1}^n \log(f_j'(f_{j,n}(1)))\right]'\\&=\sum_{j=1}^n \frac{f_j''(f_{j,n}(1))}{f_{j}'(f_{j,n}(1)) }\prod_{i=j+1}^n f_{i}'(f_{i,n}(1))=\sum_{j=1}^n \frac{f_j''(f_{j,n}(1))}{f_{j}'(f_{j,n}(1)) }\cdot\frac{f_{0,n}'(1)}{\mu_{j,n}}.
\end{align*}


\bigskip

\begin{Prf}[Theorem \ref{thm:Z-infty}]
We use the coupling from formula \eqref{X}. By Lindvall's theorem (see \cite{Lindvall}) we have that $\widetilde Z_n$ converges a.s. to a random variable $\widetilde Z_\infty$ with values in $\mathbb N_0\cup \{\infty\}$. Thus from \eqref{X} it follows that $Z_n$ has the a.s. limit
\begin{align*}
Z_\infty = \begin{cases} \widetilde Z_\infty , & \text{if } \tau_\Delta=\infty, \\
\Delta, &\text{if } \tau_\Delta < \infty. \end{cases}
\end{align*}
Moreover, we have
\begin{align*}
P[Z_\infty \neq \Delta \mid \widetilde Z_n, n \in \mathbb N_0] = P[\widetilde Z_n= Z_n \text{ for all }n \in \mathbb N_0 \mid \widetilde Z_n, n \in \mathbb N_0]\notag  = \prod_{n=1}^\infty f_n(1)^{\widetilde Z_{n-1}} ,
\end{align*}
and $\prod_{n=1}^\infty f_n(1)^{\widetilde Z_{n-1}}=0$ a.s. on the event $\{0< \widetilde{Z}_\infty < \infty\}$ if and only if $\prod_{n=1}^\infty f_n(1)=0$. Hence
\begin{align*}
P[0< Z_\infty < \infty]&= P[Z_\infty \neq \Delta, 0< \widetilde Z_\infty < \infty]= E[P[Z_\infty\neq\Delta|\widetilde{Z}_n,\ n\in\N_0];0< \widetilde Z_\infty < \infty]>0
\end{align*}
if and only if $\prod_{n=1}^\infty f_n(1)>0$ and $P[0< \widetilde{Z}_\infty < \infty]>0$. On the one hand, it is well known that $\prod_{n=1}^\infty f_n(1)>0$ if and only if
\begin{align}\label{eq:condprobab-1}
\sum_{n=1}^\infty (1-f_n(1))=\sum_{n=1}^\infty \delta_n < \infty.
\end{align}
On the other hand,  in view of Lindvall's theorem $P[0< \widetilde Z_\infty < \infty]>0$ if and only 
\begin{align}\label{eq:condprobab-2}
\sum_{n=1}^\infty (1-g_n[1]) = \sum_{n=1}^\infty \frac {f_n(1)-f_n[1]}{f_n(1)} < \infty ,
\end{align}
where recall that $g_n(s)=f_n(s)/f_n(1)$. Now, \eqref{eq:condprobab-1} implies $f_n(1) \to 1$, and therefore we obtain that \eqref{eq:condprobab-1} and \eqref{eq:condprobab-2} is equivalent to
\[ \sum_{n=1}^\infty (1-f_n[1])=\sum_{n=1}^\infty \delta_n + \sum_{n=1}^\infty  (f_n(1)-f_n[1]) < \infty. \]
This is our claim.

%
\end{Prf}

\vspace*{0.5cm}

\begin{Prf}[Theorem \ref{thm:explosion}]
Throughout this proof all the limits are taken as $n\to\infty$ unless specified otherwise.

First of all, by coupling in \eqref{X}, one sees that $\{Z_n\}_{n\in\N_0}$  and $\{\widetilde{Z}_n\}_{n\in\N_0}$ visit the same states until the process $\{Z_n\}_{n\in\N_0}$ goes to the graveyard state.

Moreover, $\{\widetilde{W}_n\}_{n\in\N_0}$, with $\widetilde{W}_n=\widetilde{Z}_n\tilde{\mu}_n^{-1}$, and $\tilde{\mu}_n=\prod_{i=1}^n f_i'(1)f_i(1)^{-1}$ is a non-negative martingale and hence $\widetilde{W}_n$ converges a.s. to an integrable random variable $\widetilde{W}\geq 0$. By applying Theorem 2 in \cite{Kersting-2020} one has that under condition \eqref{eq:cond_2},  $P[\widetilde{Z}_n\to 0]=P[\widetilde{W}=0]$. From the proof of Theorem \ref{thm:Z-infty} we have
\begin{equation*}
P[Z_n=\widetilde{Z}_n\text{ for all }n\in\N_0|\widetilde{Z}_n,\ n\in\N_0]=\prod_{n=1}^\infty f_n(1)^{\widetilde{Z}_{n-1}}\text{ a.s.},
\end{equation*}
and then,
\begin{align*}
P[Z_\infty\neq\Delta]&= E\left[\prod_{n=1}^\infty f_n(1)^{\widetilde{Z}_{n-1}}\right]\\
&=E\left[\prod_{n=1}^\infty f_n(1)^{\widetilde{Z}_{n-1}};\widetilde{Z}_n\to0\right]+E\left[\prod_{n=1}^\infty f_n(1)^{\widetilde{Z}_{n-1}};\widetilde{W}>0\right]\\
&=P[Z_n\to 0]+E\left[\prod_{n=1}^\infty f_n(1)^{\widetilde{Z}_{n-1}}\ \big|\ \widetilde{W}>0\right]P[\widetilde{W}>0],
\end{align*}
As a result,
\begin{align*}
P[\tau_a=\infty]=P[Z_\infty\neq\Delta]-P[Z_\infty=0]=E\left[\exp\Big(\sum_{n=1}^\infty\widetilde{Z}_{n-1}\log f_n(1)\Big)\ \big|\ \widetilde{W}>0\right]P[\widetilde{W}>0].
\end{align*}
Thus to have extinction or absorption in $\Delta$ with positive probability, both right-hand terms have to be positive.


First let us consider  the right-hand conditional expectation. Rewriting it as
\[E \left[\exp\Big(\sum_{n=1}^\infty\frac{\widetilde{Z}_{n-1}}{\widetilde \mu_{n-1}}\widetilde \mu_{n-1}\log f_n(1)\Big) \ \Big|\ \widetilde{W}>0\right]\]
and remembering that $\widetilde{Z}_{n}/\widetilde \mu_{n}$ converges to $\widetilde W$ a.s. we see that this  expectation is positive if and only if
\[ \sum_{n=1}^\infty\widetilde \mu_{n-1} \log f_n(1) >- \infty.\]

Second,  again due to $P[\widetilde W >0]=1-P[\widetilde Z_n \to 0]$, from Theorem 1 (iv) in \cite{Kersting-2020} it follows that $P[\widetilde W>0]>0$ if and only if
\[ \sum_{n=1}^\infty \frac{g_n''(1)}{g_n'(1)^2 \widetilde \mu_{n-1}} < \infty, \quad \text{and}\quad  \limsup_{n \to \infty}\widetilde \mu_{n}> 0,\]
which from Theorem 1 (viii) in \cite{Kersting-2020} is also equivalent to
\[ \sum_{n=1}^\infty \frac{g_n''(1)}{g_n'(1)^2 \widetilde \mu_{n-1}} < \infty \quad \text{and}\quad  \widetilde \mu_{n} \to r,\]
for some $0<r\le \infty$. This implies that $P[\widetilde W >0]>0$ is as well equivalent to the condition
\[ \sum_{n=1}^\infty \frac{g_n''(1)}{g_n'(1) \widetilde \mu_{n}} < \infty \quad \text{and}\quad  \inf_{n \ge 1}\widetilde \mu_{n}>0\]
convenient for our purpose.

Altogether we see that $P[\tau_a=\infty]>0$ holds if and only if
\[\inf_{n \ge 1}\widetilde \mu_{n}>0, \quad \sum_{n=1}^\infty \widetilde \mu_{n-1} \log f_n(1) >-\infty \quad \text{and} \quad \sum_{n=1}^\infty \frac{g_n''(1)}{g_n'(1)\widetilde \mu_{n} } < \infty.\]
Concluding the proof we rephrase this condition. First the two left-hand conditions imply that $\log f_n(1) \to 0$. Therefore we have $-\log f_n(1) \sim 1-f_n(1)=\delta_n$ (where $a_n\sim b_n$ means that there exist $c,C>0$ such that $ca_n\leq b_n\leq Ca_n$), leading to the equivalent set of conditions
\[\inf_{n \ge 1}\widetilde \mu_{n} >0, \quad \sum_{n=1}^\infty \delta_n\widetilde \mu_{n-1} < \infty \quad \text{and} \quad \sum_{n=1}^\infty \frac{g_n''(1)}{g_n'(1)\widetilde \mu_{n}} < \infty.\]
Second, these conditions entail $\sum_{n=1}^\infty \delta_n<\infty$ which means that  $\prod_{i=1}^n (1-\delta_n)= \prod_{i=1}^n f_i(1) $ converges to a positive limit. Because of $\widetilde \mu_{n}= \mu_n/\prod_{i=1}^n f_i(1)$ we end up with the conditions
\[\inf_{n \ge 1}\mu_n>0, \quad \sum_{n=1}^\infty \delta_n \mu_{n-1}< \infty \quad \text{and} \quad \sum_{n=1}^\infty \frac{f_n''(1)}{f_n'(1) \mu_n} < \infty,\]
which in turn by the same  reasoning entail the former sets of conditions. This finishes the proof.

\end{Prf}


\vspace*{0.5cm}

\begin{Prf}[Proposition \ref{thm:explosion-2}]
The proof of the first inequality is as follows. For each $k \le n$ we have
$$ P[\tau_a >n] \le P[Z_k\neq 0, Z_k \neq \Delta] \le E[Z_k] = \prod_{j=1}^k f_j'(f_{j,k}(1)) \le \mu_k .$$

For the second part, for each $t\in (0,1)$ and $n\in\N$, let us denote as $\xi_{nt}$ the random variable  with the distribution $g_{nt}(u)=f_n(ut)/f_n(t)$, $u\in [0,1]$. Then, under condition \eqref{cond-Xni}, it satisfies the following condition: there exists $c>0$ such that
\begin{equation}\label{eq:cond-t}
E[\xi_{nt}^2;\xi_{nt}\geq 2]\leq c E[\xi_{nt};\xi_{nt}\geq 2]E[\xi_{nt}|\xi_{nt}\geq 1].
\end{equation}

Indeed, if we consider the probability measures $h_t$ with weights
\begin{equation}
h_t[0]=h_t[1]=0,\quad h_t[k]=\frac{k f_n[k] t^k}{\sum_{j=2}^\infty j f_n[j]t^j},\quad k\geq 2,
\label{hmeasure}
\end{equation}
and apply Lemma 3 in \cite{Kersting-2020} for $h_{t_1}$ and $h_{t_2}$, with $t_1<t_2$, we obtain that
\begin{equation*}
\sum_{k=2}^\infty k h_t[k]=\frac{E[\xi_{nt}^2;\xi_{nt}\geq 2]}{E[\xi_{nt};\xi_{nt}\geq 2]}
\end{equation*}
is an increasing function in $t$. Thus,
\begin{equation*}
E[\xi_{nt}^2;\xi_{nt}\geq 2]\leq c E[\xi_{nt};\xi_{nt}\geq 2]\leq c E[\xi_{nt};\xi_{nt}\geq 2]E[\xi_{nt}|\xi_{nt}\geq 1],
\end{equation*}
with the last factor being greater than one.

Now, from \eqref{eq:cond-t} and Lemma 1 in \cite{Kersting-2020}, one has for each $t\in (0,1)$ and $n\in\N$
\begin{equation*}
\frac{1}{1-g_{nt}(u)}=\frac{1}{g_{nt}'(1)(1-u)}+\varphi_{nt}(u),\quad u\in [0,1],
\end{equation*}
for some function $\varphi_{nt}(\cdot)$ satisfying
\begin{equation*}
\frac{f_n''(t)f_n(t)}{2cf_n'(t)^2}\leq\frac{\varphi_{nt}(0)}{2}\leq \varphi_{nt}(u)\leq 2\varphi_{nt}(1)=\frac{f_n''(t)f_n(t)}{f_n'(t)^2},\quad u\in [0,1],
\end{equation*}
and consequently, for $0<s<t<1$,
\begin{equation}\label{eq:iter-gn-fn}
\frac{1}{f_n(t)-f_n(s)}=\frac{1}{f_n'(t)(t-s)}+\frac{\varphi_{nt}\left(s/t\right)}{f_n(t)},
\end{equation}
with
\begin{equation}\label{eq:bounds-fn}
\frac{f_n''(t)}{2cf_n'(t)^2}\leq \frac{\varphi_{nt}\left(s/t\right)}{f_n(t)}\leq\frac{f_n''(t)}{f_n'(t)^2}.
\end{equation}
By iterating the formula \eqref{eq:iter-gn-fn}, one obtains
\begin{equation*}
\frac{1}{P[\tau_a>n]}=\frac{1}{f_{0,n}(1)-f_{0,n}(0)}=\frac{1}{\prod_{j=1}^n f_j'(f_{j,n}(1))}+\sum_{j=1}^n \frac{\varphi_{j f_{j,n}(1)}\left(f_{j,n}(0)/f_{j,n}(1)\right)}{f_{j-1,n}(1)\prod_{i=1}^{j-1} f_i'(f_{i,n}(1))},
\end{equation*}
and from \eqref{eq:bounds-fn},
\begin{align}\label{eq:lower-bound-T-absorption}
\frac{1}{\prod_{j=1}^n f_j'(f_{j,n}(1))}&+\frac{1}{2c}\sum_{j=1}^n \frac{f_{j}''(f_{j,n}(1))}{f_j'(f_{j,n}(1))\prod_{i=1}^{j} f_i'(f_{i,n}(1))}\leq\nonumber\\
&\leq\frac{1}{P[\tau_a>n]}\leq \frac{1}{\prod_{j=1}^n f_j'(f_{j,n}(1))}+\sum_{j=1}^n \frac{f_{j}''(f_{j,n}(1))}{f_j'(f_{j,n}(1))\prod_{i=1}^{j} f_i'(f_{i,n}(1))}.
\end{align}
Now, the proof finishes by using the formulas in \eqref{eq:moments1} and \eqref{eq:moments2}.
\end{Prf}

\begin{Prf}[Theorem \ref{cor:asympt-E-infty}]
First note that by the Cauchy-Schwartz inequality and by \eqref{cond-Xni} one has
 \[  E[ X_{n1} ;  X_{n1} \ge 2]^2 \le  E[ X_{n1}^2 ;  X_{n1} \ge 2] \le c  E[  X_{n1} ;   X_{n1} \ge 2], \]
and therefore $ E[ X_{n1} ;  X_{n1} \ge 2] \le c$. Consequently
\begin{align}f_n'(1) \le c+1 \quad \text{and} \quad f_n''(1) \le \  E[ X_{n1}^2 ;  X_{n1} \ge 2] \le c \  E[ X_{n1} ;  X_{n1} \ge 2] \le c^2.
\label{fnprime}
\end{align}
Also, by means of the probability measure $h_t$ from \eqref{hmeasure} (with $n$ replaced with $i$ and $t=f_{i,n}(1)$)
\begin{equation}
\frac{f_i''(t)}{f_i'(t)} \le \frac{ \sum_{k=2}^\infty k^2 t^{k-2}f_i[k]}{\sum_{k=2}^\infty k t^{k-1}f_i[k]}= \frac 1t\sum_{k=2}^\infty k h_t[k]\le  \frac 1t\sum_{k=2}^\infty k h_1[k]= \frac 1t\frac{\  E[X_{i1}^2;  X_{i1}\ge 2]}{\  E[X_{i1};  X_{i1}\ge 2]} \le \frac ct ,
\label{hmeasure_2}
\end{equation}
where we again used that $\sum_{k=2}^\infty k h_t[k]$ is increasing in $t$.

Coming to our first claim note that under the condition $f_n( \rho)\ge  \rho$ for $n \ge n_0$  we have for $n_0\le i \le n$
\begin{align} f_{i,n}(1) \ge f_{i,n}( \rho) \ge  \rho.
\label{theta1}
\end{align}
The proof proceeds via induction on $i$. For $i=n$ this follows from $f_{n,n}(s)=s$, and if the claim holds true for $i>n_0$, then
\begin{align*} f_{i-1,n}(1) \ge f_{i-1,n}( \rho)= f_i(f_{i,n}( \rho)) \ge f_i( \rho) \ge  \rho.
\end{align*}
Additionally, for $0\le i \le n_0$
\begin{align} f_{i,n}(1)= f_{i,n_0}(f_{n_0,n}(1)) \ge f_{i,n_0}( \rho) >0 .
\label{theta3}
\end{align}

For $n_0\le j \le n$ we obtain
\[ \  \mu_{j,n} = \prod_{i=1}^{j} f_i'(f_{i,n}(1))  \ge \prod_{i=1}^{n_0}f_i'(f_{i,n_0}( \rho)) \prod_{i=n_0+1}^{j} f_i'( \rho). \]
This implies that there is a constant $c_1>0$ such that
\[ \  \mu_{j,n} \ge c_1 \mu_j( \rho) . \]
(For $j=n$ this is an estimate for $E[Z_n]$). With $c_1$ sufficiently small this holds for $n \le n_0$, too.
Using this estimate together with Proposition \ref{thm:explosion-2} and the formulas  \eqref{eq:lower-bound-T-absorption}, \eqref{hmeasure_2}, \eqref{theta1},  \eqref{theta3} yields
\begin{align*} \frac 1{\  P[\tau_a>n]} &\le \frac1{\  E[Z_n]} +\sum_{j=1}^{n} \frac{f_j''(f_{j,n}(1))}{f_j'(f_{j,n}(1))\   \mu_{j,n}}\\
&\le
\frac1{\  E[Z_n]} +\sum_{j=1}^{n} \frac{c}{f_{j,n}(1)\   \mu_{j,n}}\\
&\le
\frac1{c_1\mu_n( \rho)} +\sum_{j=1}^{n_0} \frac{c}{f_{j,n_0}(\rho) c_1\mu_j( \rho)}+\sum_{j=n_0+1}^{n} \frac{c}{\rho c_1\mu_j( \rho)} \\
&\le  \nu_n(\rho) \Big( \frac 1{c_1}+ \sum_{j=1}^{n_0} \frac{c}{f_{j,n_0}(\rho) c_1}+ \frac{c}{\rho c_1}\Big).
\end{align*}
 These   estimates give the first part of the theorem.

As to the second part we  notice that by our assumptions and by convexity  the generating functions $f_n(s)$ take   their values for $n\ge n_0$ and $s \ge  \sigma$ below the straight line connecting the points $( \sigma, \sigma)$ and $(1, \zeta)$ in the plane, with $\zeta= \sup_{n \ge n_0} f_n(1)$. Since by assumption $\zeta < 1$, it follows that for any $\varepsilon >0$ there is a natural number $m=m_\varepsilon$ such that
$ f_n(s) \le s-m^{-1}$ for $s \ge  \sigma+\varepsilon$.  Also, $f_n(s)\le  \sigma+ \varepsilon$ for $s \le  \sigma + \varepsilon$. Therefore, iterating these estimates, if $n \ge n_0$ and $n-i \ge m$, then
\[ f_{i,n}(1)  \le  \sigma+ \varepsilon .\]
This implies for $n\ge j > m+n_0$
\[ \  \mu_{j,n} = \prod_{i=1}^{j} f_i'(f_{i,n}(1)) \le\prod_{i=1}^{n_0} f_i'(1) \prod_{i=n_0+1}^{j-m} f_i'( \sigma+\varepsilon)  \prod_{i=j-m+1}^{j} f_i'(1).  \]
Furthermore, from \eqref{hmeasure_2} we have $f_i''( \rho) \le c f_i'( \rho)/\rho $, hence
by assumption $\inf_{i\ge 1} f_i''( \rho)>0$ also $\inf_{i\ge 1} f_i'( \rho)>0$ and consequently $\inf_{i\ge 1} f_i'( \sigma+\varepsilon)>0$.
Taking also \eqref{fnprime} into account it follows that there is a constant $c_2>0$  such that
\begin{align} 
\mu_{j,n}\le c_2  \mu_j( \sigma+\varepsilon) .
\label{exp_estimate}
\end{align}
Again, by further enlarging $c_2$, this estimate is valid for all $n \ge 1$.
Next, from
\eqref{eq:lower-bound-T-absorption} and \eqref{theta1},  for $n\geq n_0$
\[\frac 1{\  P[\tau_a>n]} \ge \frac1{\  E[Z_n]} + \frac 1{2c}\sum_{j=n_0}^n\frac{f_j''( \rho)}{f_j'(1)\  \mu_{j,n}}
 \]
 Using the assumption  $\inf_{n\ge 1} f''( \rho)>0$, \eqref{fnprime} and \eqref{exp_estimate}  we obtain
 \[ \frac 1{\  P[\tau_a>n]} \ge c_3 \sum_{j=n_0}^n \frac 1{\mu_j( \sigma+\varepsilon)} \]
 for some constant $c_3>0$. This estimate implies the assertion of the theorem, and the proof is finished.
 \end{Prf}

\vspace*{0.15cm}

\vspace*{0.5cm}

\begin{Prf}[Lemma \ref{lem:construction-branching-tree}]
We follow the steps in the proof of Lemma 1.1 in \cite{Kerting-Vatutin-2017}, thereby correcting a notation error of that proof.

First, we check that $\{g_n[d,c]\}_{1\leq d\leq c; c\in\N}$ is indeed a non-defective (proper) probability distribution.
\begin{align*}
\sum_{1\leq d\leq c<\infty} g_n[d,c]&=\frac{1}{f_{0,n}(1)-f_{0,n}(0)}\left[\sum_{c=1}^\infty f_1[c]f_{1,n}(1)^{c}-\sum_{c=1}^\infty f_1[c]f_{1,n}(0)^{c}\right]\\
&=\frac{f_1(f_{1,n}(1))-f_1(f_{1,n}(0))}{f_{0,n}(1)-f_{0,n}(0)}=1.
\end{align*}

Let $t$ be a tree of height at least $n$ and $c=c(\varnothing)$; we divide the tree $t$ into $c$ subtrees $t_1,\ldots,t_c$. Analogously, we divide the defective branching tree $T$ into subtrees $T_1,\ldots,T_c$. Thus, if $h'=h-1$, then
$$\{T\stackrel h=t\}=\{C=c,T_1\stackrel{h'}= t_1,\ldots,T_c\stackrel{h'}= t_c\},$$
and bearing in mind the independence in the reproduction among individuals
$$P_v\left[T\stackrel h= t\right]=f_1[c]\prod_{j=1}^c P_{v_1}\left[T\stackrel{h'}= t_j\right],$$
since $P_v\left[T_j\stackrel{h'}= t_j\right]=P_{v_1}\left[T\stackrel{h'}= t_j\right]$, for all $j=1,\ldots,c$. Now, note that if the tree $t$ has height at least $n$, then there is a distinguished individual $d$ in generation one, $d\in\{1,\ldots,c\}$, such that $t_d$ has height at least $n-1$, but the trees $t_j$ have heights less than $n-1$ for $j\in\{1,\ldots,d-1\}$, and the trees $t_j$, for $j\in\{d+1,\ldots,c\}$ may have height less than, equal to or greater than $n-1$. As a result 
\begin{align*}
P_v\left[T \stackrel h= t\right]=f_1[c]\prod_{j=1}^{d-1} P_{v_1}\left[T \stackrel{h'}= t_j\right] P_{v_1}\left[T\stackrel{h'}= t_d\right]\prod_{j=d+1}^c P_{v_1}\left[T\stackrel{h'}= t_j\right].
\end{align*}
For  $h \ge n$ this equation may be rewritten as 
\begin{align*}
P_v\left[T \stackrel h= t, \tau_a>n\right]&=f_1[c]\prod_{j=1}^{d-1} P_{v_1}\left[T \stackrel{h'}= t_j,\tau_0\le n-1 \right] \\&\qquad \times P_{v_1}\left[T\stackrel{h'}= t_d, \tau_a >n-1\right]\prod_{j=d+1}^c P_{v_1}\left[T\stackrel{h'}= t_j, \tau_\Delta>n-1\right].
\end{align*}
or equivalently
\begin{align*}
P_v\left[T\stackrel h=t\mid \tau_a>n\right]&=g_n[d,c] \prod_{j=1}^{d-1} P_{v_1}\left[T\stackrel{h'}= t_j \mid \tau_0\leq n-1\right]\\
&\qquad \times P_{v_1}\left[T\stackrel{h'}= t_d \mid \tau_a>n-1\right]\prod_{j=d+1}^c P_{v_1}\left[T\stackrel{h'}=t_j\mid \tau_\Delta>n-1\right].
\end{align*}
This formula implies the lemma's assertion. 
\end{Prf}


\vspace*{0.5cm}

\vspace*{0.5cm}

\begin{Prf}[Theorem \ref{thm:cond-expectation}]
First, for each $n\in\N$, by the construction we have
\begin{align*}
E[Z_n|\tau_a>n]&=1+\sum_{l=1}^{n} E[C_l-D_l]\cdot E_{v_l}[Z_{n-l}(T_{D_l+1,l})|\tau_\Delta > n-l ]\\
&=1+\sum_{l=1}^{n} E[C_l-D_l]\cdot \frac{E_{v_l}[Z_{n-l} I_{\{\tau_\Delta> n-l\}} ]}{P_{v_l}[\tau_\Delta> n-l ]}.
\end{align*}
On the one hand,
\begin{align*}
E_{v_l}[Z_{n-l} I_{\{\tau_\Delta> n-l\}} ]&=E_{v_l}\left[\tilde{Z}_{n-l} \prod_{i=1}^{n-l} f_{i+l}(1)^{\tilde{Z}_{i-1}}\right]\\
&=E_{v_l}\left[E_{v_l}\left[\tilde{Z}_{n-l} \prod_{i=1}^{n-l} f_{i+l}(1)^{\tilde{Z}_{i-1}}\Big|\tilde{Z}_{n-l-1}\right]\right]\\
&=E_{v_l}\left[E_{v_l}\left[\tilde{Z}_{n-l}\Big|\tilde{Z}_{n-l-1}\right] \prod_{i=1}^{n-l} f_{i+l}(1)^{\tilde{Z}_{i-1}}\right]\\
&=E_{v_l}\left[\tilde{Z}_{n-l-1}\frac{f_n'(1)}{f_n(1)}\prod_{i=1}^{n-l} f_{i+l}(1)^{\tilde{Z}_{i-1}}\right]\\
&\leq \prod_{i=1}^{n-l-1} f_{i+l}(1) f_n'(1) E_{v_l}\left[\tilde{Z}_{n-l-1}f_{n}(1)^{\tilde{Z}_{n-l-1}-1}\right]\\
&\leq \beta^{n-l-1} f_n'(1) E_{v_l}\left[\tilde{Z}_{n-l-1}f_{n}(1)^{\tilde{Z}_{n-l-1}-1}\right],
\end{align*}
and since the function $x\mapsto x f_n(1)^{x-1}$ has a maximum at $-1/\log(f_n(1))$, then
\begin{align*}
E_{v_l}\left[\tilde{Z}_{n-l-1}f_{n}(1)^{\tilde{Z}_{n-l-1}-1}\right]\leq \frac{1}{f_n(1)^{1+\frac{1}{\log(f_n(1))}}\log\left(\frac{1}{f_{n}(1)}\right)}.
\end{align*}

On the other hand,
\begin{align*}
P_{v_l}[\tau_\Delta> n-l]\geq P_{v_l}[\tau_\Delta=\infty]\geq P_{v_l}[\tau_0=1]=f_{l+1}(0)\geq \alpha>0.
\end{align*}

Now, we use the inequalities $E[C_l-D_l] \le E[C_l]$, $l=1,\ldots,n$. To compute the last expectation, we first determine the distribution of $C_l$, which is given by
\begin{align*}
P[C_l=c]= \sum_{d=1}^c g_{l,n}[d,c]&= \frac{f_{l,n}(1)-f_{l,n}(0)}{f_{l-1,n}(1)-f_{l-1,n}(0)}f_l[c] \sum_{d=1}^c f_{l,n}(0)^{d-1}f_{l,n}(1)^{c-d}\\
&= f_l[c]\frac{f_{l,n}(1)^c-f_{l,n}(0)^c}{f_{l-1,n}(1)-f_{l-1,n}(0)},\quad c\in\N.
\end{align*}
Therefore, for each $l=1,\ldots,n$, and $c\in\N$, using the mean value theorem, we have that for some $\xi$ between $f_{l,n}(0)^c$ and $f_{l,n}(1)^c$
\begin{align*}
\frac{P[C_l=c+k]}{P[C_l=c]} &= \frac{f_l[c+k]}{f_l[c]}\cdot \frac {f_{l,n}(1)^{c+k}-f_{l,n}(0)^{c+k}}{f_{l,n}(1)^c-f_{l,n}(0)^c}\\
& = \frac{f_l[c+k]}{f_l[c]}\cdot \frac {\big(f_{l,n}(1)^{c}\big)^{1+k/c}-\big(f_{l,n}(0)^{c}\big)^{1+k/c}}{f_{l,n}(1)^c-f_{l,n}(0)^c}\\
&= \frac{f_l[c+k]}{f_l[c]}\cdot \Big(1+ \frac kc\Big) \xi^{k/c}\\
&\le \frac{(c+k)f_l[c+k]/f_l'(1)}{cf_l[c]/f_l'(1)}.
\end{align*}

Since $\{cf_l[c]/f_l'(1)\}_{c\in\N}$ are the weights of a probability distribution, we apply \cite[Lemma 3]{Kersting-2020} and obtain
\[ E[C_l] = \sum_{c=1}^\infty c P[C_l=c] \le \sum_{c=1}^\infty c^2 \frac{f_l[c]}{f_l'(1)}= 1+\sum_{c=1}^\infty c(c-1) \frac{f_l[c]}{f_l'(1)} = 1+ \frac{f_l''(1)}{f_l'(1)} . \]

Finally, combining all the above we get
\begin{align*}
E[Z_n|\tau_a>n]&\leq1+\frac{f_n'(1)}{\alpha f_n(1)^{1+\frac{1}{\log(f_n(1))}}\log\left(\frac{1}{f_{n}(1)}\right)}\sum_{l=1}^{n} \left(1+ \frac{f_l''(1)}{f_l'(1)}\right) \beta^{n-l-1}.
\end{align*}
Using $0<\alpha \le f_n(1) \le \beta<1$ and $f_n(1)^{1/\log f_n(1)}=e$ (the Euler constant) we see that the right-hand denominator  is bounded away from 0, and our  first claim follows with
\[ c= \frac 1{e\alpha^2 \beta \log \beta^{-1}}.\]

 Finally, by \eqref{fnprime} and with the same arguments as those in the proof of \eqref{hmeasure_2} we obtain the second part.
\end{Prf}

\vspace*{0.5cm}

\end{appendix}


\section*{Acknowledgements} 

This manuscript was partially prepared while Carmen Minuesa was visiting the Institute of Mathematics, Goethe University Frankfurt, in Frankfurt am Main, and she is grateful for the hospitality and collaboration.

The authors would like to thank Serik Sagitov (Chalmers University of Technology and University of Gothenburg) for suggesting the topic of this research.

The authors also thank the anonymous referees for their valuable comments.

\vspace*{0.15cm}

\section*{Funding} 

Carmen Minuesa's research has been supported by the Ministerio de Econom\'ia y Competitividad (grant MTM2015-70522-P), Ministerio de Ciencia e Innovaci\'on (grant PID2019-108211GB-I00), the Junta de Extremadura and the European Regional
Development Fund (grants IB16099 and GR18103).


\vspace*{0.15cm}


\begin{thebibliography}{10}

\bibitem{Agresti}
A.~Agresti.
\newblock On the extinction times of varying and random environment branching
  processes.
\newblock {\em Journal of Applied Probability}, 12(1):39--46, 1975.

\bibitem{Braunsteins-Hautphenne-2019}
P.~Braunsteins and S.~Hautphenne.
\newblock Extinction in lower {H}essenberg branching processes with countably
  many types.
\newblock {\em Annals of Applied Probability}, 29(5):2782--2818, 10 2019.

\bibitem{SouzaBigginsa}
J.~C. {D'Souza} and J.~D. Biggins.
\newblock The supercritical {G}alton--{W}atson processes in varying
  environment.
\newblock {\em Stochastic Processes and their Applications}, 42:39--47, 1992.

\bibitem{Foster-Goettge-1976}
J.~H. Foster and R.~T. Goettge.
\newblock The rates of growth of the {G}alton-{W}atson process in varying
  environment.
\newblock {\em Journal of Applied Probability}, 13(1):144â€“147, 1976.

\bibitem{Fujimagari}
T.~Fujimagari.
\newblock On the extinction time distribution of a branching process in varying
  environments.
\newblock {\em Advances in Applied Probability}, 12(2):350--366, 1980.

\bibitem{Geiger-1999}
J.~Geiger.
\newblock {Elementary new proofs of classical limit theorems for
  {G}alton-{W}atson processes}.
\newblock {\em Journal of Applied Probability}, 36(2):301--309, 1999.

\bibitem{Goettge-1976}
R.~T. Goettge.
\newblock Limit theorems for the supercritical {G}alton-{W}atson process in
  varying environments.
\newblock {\em Mathematical Biosciences}, 28(1):171 -- 190, 1976.

\bibitem{Hsu-Robbins-1947}
P.~L. Hsu and H.~Robbins.
\newblock Complete convergence and the law of large numbers.
\newblock In {\em Proceedings of the National Academy of Sciences of the United
  States of America}, volume~33, pages 25--31, 1947.

\bibitem{Hu-Hu-Yin-2011}
Y.~Hu, W.~Hu, and Y.~Yin.
\newblock The extinction of a branching process in a varying or random
  environment.
\newblock In S.~Li, X.~Wang, Y.~Okazaki, J.~Kawabe, T.~Murofushi, and L.~Guan,
  editors, {\em Nonlinear Mathematics for Uncertainty and its Applications.
  Advances in Intelligent and Soft Computing}, volume 100, pages 309--315.
  Springer, Berlin, Heidelberg, 2011.

\bibitem{Ispany-2016}
M.~Isp\'any.
\newblock Some asymptotic results for strongly critical branching processes
  with immigration in varying environment.
\newblock In I.~del Puerto, M.~Gonz\'{a}lez, C.~Guti\'{e}rrez, R.~Mart\'{i}nez,
  C.~Minuesa, M.~Molina, M.~Mota, and A.~Ramos, editors, {\em Branching
  Processes and Their Applications}, volume 219 of {\em Lecture Notes in
  Statistics}, pages 77--95. Springer, 2016.

\bibitem{Iwasa-Michor-Nowak-2004}
Y.~Iwasa, F.~Michor, and M.~A. Nowak.
\newblock Evolutionary dynamics of invasion and escape.
\newblock {\em Journal of Theoretical Biology}, 226(2):205--214, 2004.

\bibitem{Jagers-1974}
P.~Jagers.
\newblock {G}alton-{W}atson processes in varying environments.
\newblock {\em Journal of Applied Probability}, 11(1):174--178, 1974.

\bibitem{Karlin-Tavare-1982}
S.~Karlin and S.~TavarÃ©.
\newblock Linear birth and death processes with killing.
\newblock {\em Journal of Applied Probability}, 19(3):477--487, 1982.

\bibitem{Keiding-Nielsen-1975}
N.~Keiding and J.~E. Nielsen.
\newblock Branching processes with varying and random geometric offspring
  distributions.
\newblock {\em Journal of Applied Probability}, 12(1):135â€“141, 1975.

\bibitem{Kersting-2020}
G.~Kersting.
\newblock A unifying approach to branching processes in a varying environment.
\newblock {\em Journal of Applied Probability}, 57(1):196--220, 2020.

\bibitem{Kerting-Vatutin-2017}
G.~Kersting and V.~Vatutin.
\newblock {\em Discrete Time Branching Processes in Random Environment}.
\newblock ISTE Ltd and John Wiley and Sons, Inc., 2017.

\bibitem{Lindvall}
T.~Lindvall.
\newblock Almost sure convergence of branching processes in varying and random
  environments.
\newblock {\em The Annals of Probability}, 2(2):344--346, 1974.

\bibitem{MaCphee-Schuh-1983}
I.~M. MacPhee and H.~J. Schuh.
\newblock A {G}alton-{W}atson branching process in varying environments with
  essentially constant offspring means and two rates of growth.
\newblock {\em Australian Journal of Statistics}, 25(2):329--338, 1983.

\bibitem{Robertson-1978}
A.~Robertson.
\newblock The time of detection of recessive visible genes in small
  populations.
\newblock {\em Genetics Research}, 31:255--264, 1978.

\bibitem{Sagitov-Minuesa-2017}
S.~Sagitov and C.~Minuesa.
\newblock Defective {G}alton--{W}atson processes.
\newblock {\em Stochastic Models}, 33(3):451--472, 2017.

\bibitem{Yu-Pei-2009}
J.~Yu and J.~Pei.
\newblock Extinction of branching processes in varying environments.
\newblock {\em Statistics and Probability Letters}, 79(17):1872 -- 1877, 2009.

\end{thebibliography}

\end{document}